\documentclass[12pt]{article}
\begin{document}
\title{Regularities of Twin, Triplet and Multiplet Prime Numbers}
\author{H. J. Weber\\Department of Physics\\University of Virginia\\
Charlottesville, VA 22904, U.S.A.}
\maketitle
\begin{abstract}
Classifications of twin primes are established and then applied to 
triplets that generalize to all higher multiplets. Mersenne and 
Fermat twins and triplets are treated in this framework. Regular 
prime number multiplets are related to quadratic and cubic prime 
number generating polynomials. 
\end{abstract}
\vspace{3ex}
\leftline{MSC: 11N05, 11N32, 11N80}
\leftline{Keywords: Twin primes, triplets, regular multiplets.} 


\section{Introduction}

There exist extensive tables of twin, triplet and quartet 
primes. There are no such systematic analyses of higher 
multiplets. Here we outline a more systematic analysis 
of generalized twin primes, triplets and regular 
multiplets that are connected with prime producing 
quadratic and cubic polynomials. We follow standard 
practice ignoring as trivial the prime pairs $(2,p)$ 
of odd distance $p-2$ with $p$ any odd prime. 

\section{Classifications of twin and triplet primes}

{\bf Definition~2.1.} A triplet $p_i, p_m=p_i+2d_1, 
p_f=p_m+2d_2$ of odd prime numbers with $p_i<p_m<p_f$ 
is called a {\it generalized triplet.} 

Each triplet consists of three {\it generalized twin 
primes} $(p_i, p_m), (p_m, p_f),\\(p_i, p_f)$. 
Empirical laws governing triplets therefore are 
intimately tied to those of the generalized twin 
primes. 

\subsection{Twin prime classifications}

There are two schemes of parametrizations of twin primes 
which we now describe.  

{\bf Theorem~2.2.} {\it Let $2D$ be the distance between odd  
prime numbers $p_i,p_f$ of the pair, $D$ a natural number. 
Then there are three mutually exlusive classes of generalized 
twin primes that are parametrized as follows.} 
\begin{eqnarray}
I:~p_i&=&2a-D,~p_f=2a+D,~D~\rm{odd};\\
II:~p_i&=&3(2a-1)-D,~p_f=3(2a-1)+D,~2|D,~3\not|D;\\
III:~p_i&=&2a+1-D,~p_f=2a+1+D,~6|D,~D\geq 6,
\end{eqnarray}
{\it where $a$ is the running integer variable. Values of 
$a$ for which a prime pair of distance $2D$ is reached 
are unpredictable (called arithmetic chaos).}

Each of these three classes of this classification~\cite{hjw}
contains infinitely many (possibly empty) subsets of prime 
pairs at given even distance.

{\it Proof.} Let us first consider the case of odd $D.$ Then 
$p_i=2a-D$ for some positive integer $a$ and, therefore, 
$p_f=p_i+2D=2a+D.$ The median $2a$ is the running integer 
variable of this class $I$. 

For even $D$ with $D$ not divisible by $3,$ let $p_i=2n+1-D$ 
so that $p_f=2n+1+D$ for an appropriate integer $n.$ Let 
$p_i\neq 3,$ thus excluding a possible first pair with $p_i=3$ 
as {\it special.} Since of three odd natural numbers at 
distance $D$ from each other one is divisible by $3,$ the 
median $2n+1$ must be so, hence $2n+1=3(2a-1)$ for an 
appropriate integer $a.$ Thus, the median $3(2a-1)$ of the 
pair $3(2a-1)\pm D$ is again a linear function of a running 
integer variable $a$. These prime number pairs constitute 
the 2nd class $II$. 

This argument is not valid for prime number pairs 
with $6|D,$ but they can obviously be parametrized as 
$2a+1\pm 6d,~D=6d.$ They comprise the 3rd and last class 
$III$ of generalized twins. Obviously, these three 
classes are mutually exclusive and complete except for 
the special cases.~$\diamond$   
  
{\it Example 1.} 

Ordinary twins $2a\pm 1$ for $a=2, 3, 6, 9, 15, \ldots$ have 
$D=1$ and are in class $I$. For $D=3,~2a\pm 3$ are twins for 
$a=4, 5, 7, 8, 10, \ldots.$ For $D=5,~2a\pm 5$ are twins for 
$a=4, 6, 9, 12, 18, \ldots.$ No twins are ever missed or 
special in this class $I$, an advantage of this 
classification.  

For $D=2,~3(2a-1)\pm 2$ are twins in class $II$ for $a=2, 3, 
4, 7, 8, \ldots.$

For $D=6,~2a+1\pm 6$ are twins in class $III$ for $a=5, 6, 
8, 11, 12, \ldots.$ 

Special twins are $5\pm 2,~7\pm 4,~11\pm 8,\ldots$ 
 
The {\it second classification} of generalized twins 
involves arithmetic progressions of conductor $6$ as 
their regular feature. It is well known that, except 
for the first pair $3, 5,$ ordinary twins all have the 
form $(6m-1, 6m+1)$ for some natural number $m.$ They 
belong to class $I.$           
     
{\it Example 2.}

Prime pairs at distance $2D=4$ are in class $II$ and of 
the form $6m+1, 6(m+1)-1$ for $m=1, 2, 3, \dots$ except 
for the singlet $3, 7.$ At distance $2D=6$ they are in 
class $I$ and have the form $6m-1, 6(m+1)-1$ for $m=1, 
2, 3, \ldots$ that are intertwined with $6m+1, 6(m+1)+1$ 
for $m=1, 2, 5, \dots.$ At distance $2D=12$ they are in 
class $III$ and of the form $6m-1, 6(m+2)-1$ for $m=1, 
3, 5, \ldots$ intertwined with $6m+1, 6(m+2)+1$ for $m=1, 
3, 5, \ldots.$ 

In general, the rules governing the form $6m\pm 1, 6m+b$ 
depend on the arithmetic of $D$ and $a$ making the second  
classification of generalized twin primes difficult to 
deal with generally. We now apply Theor.~2.2.    

{\bf Theorem~2.3.} {\it Let $2D$ be the distance between 
odd prime numbers $p_i,p_f$ of the pair. Then for class III, 
$D\equiv 0 \pmod{6}$ and $p_{f,i}\equiv 2a+1 \pmod{6}.$ If 
$a\equiv 0 \pmod{3}$ then $p_{f,i}\equiv 1 \pmod{6}.$ If 
$a\equiv -1 \pmod{3}$ then $p_{f,i}\equiv -1 \pmod{6}.$ 

For class II, $D\equiv 2r \pmod{6},~r=\pm 1$ yields 
$p_i\equiv \pm 1 \pmod{6}, p_f\equiv \mp1\pmod{6}; r=\pm 2$ 
gives $p_i\equiv \mp 1,~p_f\equiv \pm 1\pmod{6}.$
 
For class $I,$ and $D\equiv 1+2r\pmod{6},~r=0,\pm 1,~
a\equiv a_0\pmod{3}$ the prime pair obeys 
$p_i\equiv -1-2r-2a_0 \pmod{6}, p_f\equiv 1+2r+2a_0
\pmod{6}.$ For $r=0$ and $r=-1,$ i.e. $D\equiv \pm 1
\pmod{6},~a_0=0.$ For $r=1,$ i.e.} $D\equiv 3\pmod{6},
~a\not\equiv 0\pmod{3}.$     

{\it Proof.} For class $III,~a\equiv 1 \pmod{3}$ is 
ruled out because then $p_{f,i}\equiv 3\pmod{6}.$ 
For class $II,$ $p_{f,i}=3(2a-1)\pm D,$ so 
$D\equiv 2r\pmod{6}$ implies $p_{f,i}\equiv -3\pm 
2r\pmod{6}$ for $r=\pm 1,$ etc. For class $I,$ and 
$D\equiv 1\pmod{6}$ the cases $a\equiv \pm 1\pmod{3}$ 
are ruled out because they imply either $3|p_i$ or 
$3|p_f,$ except $a=2,$ i.e. $D=2.$ This also is the 
case for $D\equiv -1\pmod{6}.$ For $D\equiv 3\pmod{6}$ 
the case $a\equiv 0\pmod{3}$ is obviously ruled 
out.~$\diamond$  
     
This concludes the classifications of twin primes  
at constant distances.  

\subsection{Triplet prime classifications}

Generalized triplets have the form $6m\pm 1, 6m+a_1, 
6m+a_2,$ except for singlet exceptions and this 
generalizes to prime quadruplets, quintuplets, etc. 

{\it Rules} for {\it singlets} or {\it exceptions} among 
generalized triplet primes are the following. 

{\bf Theorem~2.4.} {\it (i) There is at most one 
generalized prime number triplet with distances 
$[2D, 2D]$ for $D=1, 2, 4, 5, \ldots$ and $3\not| D.$ 

(ii) When the distances are $[2d_1, 2d_2]$ with 
$3|d_2-d_1,$ and $3\not| d_1,$ there is at most one 
triplet $p_i=3, p_m=3+2d_1, p_f=3+2d_1+2d_2$ for 
appropriate integers} $d_1, d_2$.     

{\it Proof.} (i) Because of three odd numbers in a row 
one is divisible by $3,~3, 5, 7$ is the only triplet at 
distances $[2, 2]$ and, for the same reason, there is 
only one triplet $3, 7, 11$ at distances $[4, 4],$ one 
only at distances $[8, 8]$ i.e. $3, 11, 19$, at 
$[10, 10]$ i.e. $3, 13, 23$ and in general at distances 
$[2D, 2D]$ for $D$ not divisible by $3$. The argument 
fails when $3|D.$ (ii) Of $p_i,~p_m\equiv p_i+2d_1 
\pmod{3},~ p_f\equiv p_i+4d_1 \pmod{3}$ at least one 
is divisible by $3,$ which must be $p_i.$~$\diamond$

Naturally, the question arises: Are there infinitely  
many such singlets, i.e. exceptional triplets? 

{\it Example 3.} 

At distances $[2, 8],$ the singlet is $3, 5, 13$ and at 
$[8, 2],$ it is $3, 11, 13.$ 

At distances $[2, 4],$ the triplets are $2n-3, 2n-1, 2n+3$ 
with $n\equiv 1 \pmod{3}$ and at $[4,2],$ they are $2n-3, 
2n+1, 2n+3$ with $n\equiv -1 \pmod{3}.$ Writing $n=3m\pm 1$ 
for these cases, we obtain $6m-1, 6m+1, 6m+5$ and $6m-5, 
6m-1, 6m+1$ for these triplets, respectively. At distances 
$[2, 4],$ triplets occur for $n=4, 7, 10, \ldots$ i.e. 
$m=1, 2, 3, \ldots,$ while at $[4, 2],$ they are at 
$n=5, 8, 20, \ldots$ i.e. $m=2, 3, 7, \ldots.$ These 
triplets are in the classes $(I, II)$ and $(II, I),$ 
respectively, with $I, II$ denoting symbolically the 
classes of the first generalized twin prime classification.      

At distances $[2, 6]$ the generalized triplets are 
$2n-3, 2n-1, 2n+5$ with $n\equiv 1 \pmod{3}$ or 
$6m-1, 6m+1, 6m+7.$ The triplet $3, 5, 11$ for $n=3$ 
is the only exception. For $[6, 2]$ they are $2n-5, 
2n+1, 2n+3$ with $n\equiv -1 \pmod{3}$ or $6m-1, 6m+5, 
6m+7.$ They all are in the class $(I, I).$ 

Applying Theors.~2.2, 2.3 we obtain the following 
triplet classifications.  

{\bf Corollary~2.5.} {\it (i) The class $(I, I)$ is 
made up of the triplets 
\begin{eqnarray}
p_i=2a-D_1,~p_m=2a+D_1=2b-D_2,~p_f=2b+D_2
\end{eqnarray} 
with odd $D_1, D_2$ and $a, b$ appropriate integers 
subject to $b-a=(D_1+D_2)/2.$ Hence 
\begin{eqnarray}
p_f=2a+D_1+2D_2,
\end{eqnarray} 
and the prime number pair 
\begin{eqnarray}
(p_f, p_i)=2a+D_2\pm (D_1+D_2)
\end{eqnarray}
is in class $II,$ or special, or $III$.

(ii) If} 
\begin{eqnarray}\nonumber
&&D_1\equiv 2r_1+1 \pmod{6}, r_1=0,\pm 1;~ 
D_2\equiv 2r_2+1 \pmod{6}, r_2=0,\pm 1,\\&&
a\equiv a_0 \pmod{3}, a_0=0,\pm 1
\end{eqnarray}
{\it then} 
\begin{eqnarray}\nonumber
&&p_i\equiv -2r_1-1+2a_0 \pmod{6},~p_m\equiv 
2r_1+1+2a_0 \pmod{6},\\\nonumber &&p_f\equiv 
2a_0+3+2r_1+2r_2 \pmod{6},~D_1+D_2\equiv 
2(r_1+r_2)+2 \pmod{6}.\\
\end{eqnarray}   
Now specific cases (ii) can be worked out by 
substituting values for $r_i, a_0$: $D_1\equiv 
1\pmod{6},~a\equiv 0\pmod{3}$ yield $p_i\equiv 
-1\pmod{6},~p_m\equiv 1\pmod{6},~p_f\equiv 
1+2D_2\pmod{6}\equiv \pm 1\pmod{6}$ with 
$D_2\equiv 1\pmod{6}$ ruled out, etc. Of 
course, Example~3 is consistent with this.      
 
{\bf Corollary~2.6.} {\it (i) The class $(I, II)$ 
consists of prime triplets at distances $D_1$ 
odd and $D_2$ even so that
\begin{eqnarray} 
p_i=2a-D_1,~p_m=2a+D_1=3(2b-1)-D_2,~p_f=3(2b-1)+D_2 
\end{eqnarray}
with appropriate $a, b$ subject to $a=3b-\frac{1}{2}
(D_1+D_2+3).$ Hence $p_f=2a+D_1+2D_2$ and the twin 
\begin{eqnarray}
(p_f, p_i)=2a+D_2\pm (D_1+D_2)
\end{eqnarray}
is in class $I.$  

(ii) If}
\begin{eqnarray}\nonumber
&&D_1\equiv 2r_1+1 \pmod{6}, r_1=0,\pm 1;~ 
D_2\equiv 2r_2 \pmod{6}, r_2=\pm 1,\\
&&a\equiv a_0 \pmod{3}, a_0=0,\pm 1
\end{eqnarray}
{\it then} 
\begin{eqnarray}\nonumber
&&p_i\equiv -2r_1-1+2a_0 \pmod{6},~p_m\equiv 
2r_1+1+2a_0 \pmod{6}\\&&\equiv -3-2r_2 \pmod{6},
~p_f\equiv -3+2r_2 \pmod{6}. 
\end{eqnarray}

{\bf Corollary~2.7.} {\it The class $(II, II)$ has 
$D_1$ even and $D_2$ even with both $D_i$ not 
divisible by $3.$ The triplets are 
\begin{eqnarray}\nonumber
p_i&=&3(2a-1)-D_1,~p_m=3(2a-1)+D_1=3(2b-1)-D_2,\\
p_f&=&3(2b-1)+D_2 
\end{eqnarray}
so that $3(b-a)=\frac{1}{2}(D_1+D_2).$ Therefore 
$3|D_1+D_2$ and $D_1\equiv -D_2 \pmod{3}.$ Since 
\begin{eqnarray}\nonumber
p_f&=&p_i+2(D_1+D_2)=3(2a-1)+D_1+2D_2,\\ 
(p_f, p_i)&=&3(2a-1)+D_2\pm (D_1+D_2)
\end{eqnarray} 
is in class $III$ with $6|D_1+D_2.$ 

(ii) If} 
\begin{eqnarray}
D_i\equiv 2r_i \pmod{6},~r_i=\pm 1 
\end{eqnarray} 
{\it then} 
\begin{eqnarray}\nonumber
&&p_i\equiv -3-2r_1 \pmod{6},~p_m\equiv -3+2r_1 
\pmod{6}\equiv -3-2r_2 \pmod{6},\\&&p_f\equiv 
-3+2r_2 \pmod{6}\equiv p_i \pmod{6}. 
\end{eqnarray}
{\bf Corollary~2.8.} {\it The class $(I, III)$ 
has the triplets} 
\begin{eqnarray}
p_i=2a-D_1,~p_m=2a+D_1=2b+1-6d_2,~p_f=2b+1+6d_2 
\end{eqnarray}
{\it where $2b+1=2a+D_1+6d_2.$ So 
$p_f=2a+D_1+12d_2$ and} 
\begin{eqnarray}
(p_f, p_i)=2a+6d_2\pm (D_1+6d_2)
\end{eqnarray}
{\it is in class $I.$

(ii) If} 
\begin{eqnarray}
D_1\equiv 1+2r_1 \pmod{6},~r_1=0, \pm 1; a\equiv 
a_0 \pmod{3}
\end{eqnarray}
{\it then}
\begin{eqnarray}\nonumber
&&p_i\equiv -1-2r_1+2a_0 \pmod{6},~p_m\equiv 
2b+1 \pmod{6}\equiv\\\nonumber&&1+2r_1+2a_0
\pmod{6},~p_f\equiv p_m\pmod{6},~2b\equiv 
2(r_1+a_0)\pmod{6}.\\
\end{eqnarray} 

{\bf Corollary~2.9.} {\it Class $(II, III)$ has 
$D_1$ even and not divisible by $3$ and $D_2=6d_2$ 
with the triplet form 
\begin{eqnarray}\nonumber
p_i&=&3(2a-1)-D_1,~p_m=3(2a-1)+D_1=2b+1-6d_2,\\
p_f&=&2b+1+6d_2
\end{eqnarray}
so that 
\begin{eqnarray}
2b+1=3(2a-1)+D_1+6d_2,~p_f=3(2a-1)+D_1+12d_2.
\end{eqnarray}
So 
\begin{eqnarray}
(p_f, p_i)=3(2a-1)+6d_2\pm (D_1+6d_2)
\end{eqnarray}
is in class $II.$
 
(ii) If}
\begin{eqnarray}
D_1\equiv 2r_1 \pmod{6}, r_1=\pm 1,~a\equiv a_0 
\pmod{3}, a_0=0, \pm 1
\end{eqnarray}
{\it then}
\begin{eqnarray}\nonumber
&&p_i\equiv -3-2r_1 \pmod{6},~p_m\equiv -3+2r_1 
\pmod{6}\\&&\equiv 2b+1 \pmod{6},~p_f\equiv p_m 
\pmod{6}.   
\end{eqnarray}

{\bf Corollary~2.10.} {\it The class $(III, III)$ 
has the triplet form 
\begin{eqnarray}
p_i=2a+1-6d_1,~p_m=2a+1+6d_1=2b+1-6d_2, 2b+1+6d_2
\end{eqnarray}
for appropriate $a, b$ so that
\begin{eqnarray} 
b=a+3(d_1+d_2),~p_f=2a+1+6(d_1+2d_2)
\end{eqnarray}
and 
\begin{eqnarray}
(p_f, p_i)=2a+1+6d_2\pm 6(d_1+d_2) 
\end{eqnarray}
is in class $III,$ too. 

(ii) If $a\equiv a_0 \pmod{3},~b\equiv b_0 
\pmod{3}$ then} 
\begin{eqnarray}
p_i\equiv 2a_0+1 \pmod{6}\equiv p_m \pmod{6}
\equiv p_f \pmod{6}. 
\end{eqnarray}
The classes $(II, I), (III, I), (III, II)$ are 
handled similarly. Several examples have been 
given above. These nine classes of generalized 
prime number triplets are mutually exclusive and 
complete except for the singlets. 

These twin and triplet prime classifications 
represent regularities that generalize to  
quadruplet primes which come in $3^3$ mutually 
exlusive classes, quintuplet primes in $3^4$ 
such classes except for singlets, etc.         
    
\section{Special twin and triplet primes}

Here we consider Mersenne and Fermat twins and 
triplets. 

\subsection{Mersenne twins}
  
A simple application of the second classification 
is the following  

{\bf Corollary~3.1.} {\it If $2^p-1,$ with an odd 
prime number $p,$ is a Mersenne prime, then $2^p+1$ 
is composite.} 

{\it Proof.} Since $2^p-1\neq 6m-1,$ the pair 
$2^p\pm 1$ is not a twin prime.~$\diamond$ 

Of course, it is well known that $3|2^p+1$ but this 
requires an algebraic identity for the factorization: 
\begin{eqnarray}
a^p+1=(a+1)(a^{p-1}-a^{p-2}\pm\ldots +a^2-a+1).
\end{eqnarray}    

Let us now consider Mersenne twins with $2^p-1$ 
as the {\it first} prime number of the pair.  

{\it Example~4.}

The pair $2^p-1, 2^p+5$ is a Mersenne pair in 
class $I$ for $p=3, 5, \ldots,~p\equiv -1
\pmod{3}$ and $p\geq 7.$ The qualification is 
due to the factorization in Lemma~3.2.   

The pair $2^p-1, 2^p+9$ is a Mersenne pair in 
class $I$ for $p=2, 3, 5, 7, \ldots.$

The pair $2^p-1, 2^p+3$ is a Mersenne pair in 
class $II$ for $p=2, 3, 7, \ldots,\\p\equiv -1
\pmod{4}.$ The restriction is due to Lemma~3.3.

The pair $2^p-1, 2^p+11$ is a Mersenne pair in 
class $III$ for $p=3, 5, 7, \ldots.$ 

{\bf Lemma~3.2.} {\it If} $p\equiv 1\pmod{3}$ 
{\it and} $p\geq 7$ {\it then} $7|2^p+5.$ 

{\it Proof.} This follows from the factorization 
\begin{eqnarray}
2^p+2^3-2-1=(2^3-1)(2^{p-3}+2^{p-6}+\cdots +2+1).
~\diamond
\end{eqnarray}

{\bf Lemma~3.3.} {\it If} $p\equiv 1\pmod{4}$ 
{\it then} $5|2^p+3.$ 

{\it Proof.} 
\begin{eqnarray}
2^p+2^2-2+1=(2^2+1)\left(\sum_{j=1}^{(p-1)/2}
(-1)^{j-1}2^{p-2j}+1\right).~\diamond
\end{eqnarray}

Next we list Mersenne twins with $2^p-1$ as the 
{\it second} prime number of the pair.  

{\bf Proposition~3.4.} {\it (i) The pair} $2^p-5, 
2^p-1$ {\it is a Mersenne pair in class $II$ for 
$p=3$ only. There are no Mersenne prime twins} 
$2^p-2^{2n+1}-3, 2^p-1;~n=1, 2,\ldots$ {\it except} 
$p=3, n=0.$ {\it (ii) There are no Mersenne twins} 
$2^p-3, 2^p-1$ {\it when} $p\equiv -1\pmod{4}$ {\it 
except for the pair $5, 7.$ (iii) There are 
no Mersenne twins} $2^p-7,2^p-1$ {\it when} $p\equiv 
1\pmod{4}.$   

{\it Proof.} (i) This holds because $3|2^p-5$ for 
$p\geq 3$ which follows from the first factorization 
\begin{eqnarray}\nonumber
&&2^p-2-2-1=(2+1)\left(\sum_{m=1}^{p-1}(-2)^{p-m}
-1\right);\\\nonumber&&2^p-2^{2n+1}-2-1=(2+1)\left(
\sum_{m=1}^{p-2n-1}(-2)^{p-m}-1\right),\\&&
p\geq 2n+3,~n=1,2,\dots.
\end{eqnarray} 
and the next cases from the second factorization. 
The case $n=1$ gives $3|2^p-11,~p\geq 5.$ 

(ii) is due to the factorization, for $p\equiv -1
\pmod{4},$ 
\begin{eqnarray}
2^p-2^2+2-1=(2^2+1)\left(\sum_{m=1}^{(p-1)/2}
(-1)^{m-1}2^{p-2m}-1\right). 
\end{eqnarray}
(iii) follows from the factorization, for 
$p\equiv 1 \pmod{4},$ 
\begin{eqnarray}
2^p-2^2-2-1=(2^2+1)\left(\sum_{m=1}^{(p-1)/2}
(-1)^{m-1}2^{p-2m}-1\right).~\diamond 
\label{fac3}
\end{eqnarray}

{\it Example~5.}

The pair $2^p-15, 2^p-1$ is in class $I,$ and 
$p=5, 7$ are such cases. 

The pair $2^p-19, 2^p-1$ is in class $I,$ and 
$p=5, 7$ are such cases.

The pair $2^p-13, 2^p-1$ is in class $III,$ and 
$p=5, 13$ are such cases.

The pair $2^p-25, 2^p-1$ is in class $III,$ and 
$p=5, 7, 13$ are such cases.

\subsection{Fermat twins}

Here we consider Fermat prime pairs with the 
Fermat prime being its first member. 

{\it Example~6.}

$2^{2^n}+1, 2^{2^n}+3$ are twins in class $I$ for
$n=0, 1, 2, 4,\ldots,$ i.e. $2^{n-1}\not\equiv 1
\pmod{3},~n>1.$ The qualification is due to (i) in 
Lemma~3.6. 

The pair $2^{2^n}+1, 2^{2^n}+7$ is in class $I$ 
and $n=1, 2, 3$ are such cases. 

The pair $2^{2^n}+1, 2^{2^n}+13$ is in class $II$ 
and $n=1, 2, 3$ are such cases.  

{\bf Corollary~3.5.} {\it (i) The twin prime} 
$3,~7,$ {\it is the only one of} $2^{2^n}+1, 
2^{2^n}+5$ {\it in class $II$ for} $n=0.$ {\it (ii) 
The pair} $2^{2^n}+1, 2^{2^n}+9$ {\it is in class 
$III$ and $n=0, 1$ are the only such cases.}

{\it Proof.} (i) follows from the factorization  
\begin{eqnarray}
2^{2^n}+2+2+1=(2+1)\left(\sum_{j=1}^{2^n-1}(-1)^{j-1}
2^{2^n-j}+1\right),~n>0.
\end{eqnarray}
and (ii) from 
\begin{eqnarray}
2^{2^n}+2^2+2^2+1=(2^2+1)\left(\sum_{j=1}^{2^{n-1}-1}
(-1)^{j-1}2^{2^n-2j}+1\right),~n>0.~\diamond
\end{eqnarray}
{\bf Lemma~3.6.} {\it (i) If} $2^n\equiv 2\pmod{3},
~n>1$ {\it then} $7|2^{2^n}+3.$ {\it (ii) If} 
$2^n\equiv 1\pmod{3}$ {\it then} $7|2^{2^n}+5.$ 

{\it Proof.} (i) follows from the factorization
\begin{eqnarray}\nonumber
&&2^{2^n}+2^3-2^2-1=(2^3-1)\left(2^{2^n-3}+2^{2^n-6}+
\cdots +2^2+1\right),\\&&2^n\equiv 2\pmod{3},
\end{eqnarray}
and (ii) from 
\begin{eqnarray}\nonumber
&&2^{2^n}+2^3-2-1=(2^3-1)\left(2^{2^n-3}+2^{2^n-6}+
\cdots +2+1\right),\\&&2^n\equiv 1\pmod{3}.~\diamond 
\end{eqnarray}
Fermat prime pairs with the Fermat prime being its 
second member are the following.

{\it Example~7.}

The pair $2^{2^n}-5, 2^{2^n}+1$ is in class $I$ and 
$n=2, 3$ are such cases.

The pair $2^{2^n}-3, 2^{2^n}+1$ is in class $II$ and 
$n=2$ is such a case.

The pair $2^{2^n}-11, 2^{2^n}+1$ is in class $III$ 
and $n=2$ is such a case.

{\bf Proposition~3.7.} {\it There is no Fermat twin 
primes of the form} $2^{2^n}-1, 2^{2^n}+1; 2^{2^n}-7, 
2^{2^n}+1;~2^{2^n}-19, 2^{2^n}+1; 2^{2^n}+1;~
2^{2^n}-21$. 

{\it Proof.} This holds because $5|2^{2^n}-1, 
3|2^{2^n}-7, 2^{2^n}-19, 5|2^{2^n}-21$ which is 
based on the following factorizations: 
\begin{eqnarray}\nonumber
&&2^{2^n}-1=(2^2+1)\left(\sum_{j=1}^{2^{n-1}-1}
(-1)^{j-1}2^{2^n-2j}-1\right),~n>1,\\\nonumber&&
2^{2^n}+(-1)^{m+1}2^{2m}-2^2-1=(2^2+1)\left(
\sum_{j=1}^{2^{n-1}-m}(-1)^{j-1}2^{2^n-2j}-1\right),
\\\nonumber&& 1\leq m<2^{n-1},~n>1,\\\nonumber&&
2^{2^n}-2^{2m}-2-1=(2+1)\left(\sum_{j=1}^{2^n-2m}
(-1)^{j-1}2^{2^n-j}-1\right),\\&&m=1, 2,\ldots. 
~\diamond
\end{eqnarray}

\subsection{Mersenne triplets}

Here we list prime triplets where the Mersenne prime 
comes {\it first.}

{\it Example~8.}

$2^p-1,~2^p+3,~2^p+9$ yields  triplets for $p=2,~3: 
3,~7,~13; 7,~11,~13.$ 

$2^p-1,~2^p+5,~2^p+9$ yields triplets for $p=3,~5: 
7,~13,~17; ~31,~37,~41$ and $2^p-1,~2^p+5,~2^p+11$ 
yields triplets for $p=3,~5: 7,~13,~19; ~31,~37,~43.$

$2^p-1,~2^p+3,~2^p+11$ yields triplets for $p=3,~7: 
7,~11,~19; ~127,~131,~139.$

{\bf Corollary~3.8.} {\it (i)} $2^p-1,~2^p+3,~2^p+7$ 
{\it for} $p=2$ {\it yields the only such triplet} 
$3,~7,~11.$ {\it (ii)} $2^p-1,~2^p+5,~2^p+7$ 
{\it yields no Mersenne triplets, and (iii)} 
$2^p-1,~2^p+3,~2^p+5$ {\it which is in class $II,I$ 
yields no triplets except for} $p=3,$ {\it namely} 
$7,~11,~13.$  

{\it Proof.} (i) and (ii) follow from the 
factorization   
\begin{eqnarray}
2^p+2^2+2+1=(2+1)\left(\sum_{j=1}^{p-2}(-2)^{p-j}
+1\right),~p>2.
\label{fac4}
\end{eqnarray}
(iii) follows from Lemma~3.2 and the factorization   
\begin{eqnarray}\nonumber
2^p+2^3-2^2-1&=&(2^3-1)(2^{p-3}+2^{p-6}+\cdots +2^2+1), 
\\p&\equiv& -1\pmod{3}.~\diamond
\end{eqnarray} 

Prime triplets where the Mersenne prime comes 
{\it last} are the following.  

{\it Example~9.}

For $p=3$ the triplet $2^p-5,~2^p-3,~2^p-1$ 
becomes $3,~5,~7$ which is a singlet case. 

$2^p-7,~2^p-3,~2^p-1$ yields $23,~29,~31$ for $p=5.$ 

$2^p-13,~2^p-3,~2^p-1$ yields a triplet for $p=5: 
19,~29,~31.$ 

$2^p-19,~2^p-15,~2^p-1$ yields triplets 
for $p=5,~7: 13,~17,~31;~109,~113,~127.$ 

{\bf Corollary~3.9.} $2^p-7,~2^p-5,~2^p-1$ {\it yields 
no Mersenne triplets. $2^p-1,~2^p+3,~2^p+5$ yields no 
triplets except for $p=2.$}

{\it Proof.} This follows from the factorizations 
(31), (32), (33),~(\ref{fac3}).$~\diamond$

{\bf Corollary~3.10.} {\it The triplet $2^p-3, 2^p-1, 
2^p+3$ is the only one.}

{\it Proof.} This follows from the factorization of 
$5|2^p+3$ in Eq.~(3.3) of Lemma~3.3 for 
$p\equiv 1\pmod{4}$ and $5|2^p-3$ for $p\equiv 
-1\pmod{4}$ in Prop.~3.4.~$\diamond$  

Finally, prime triplets where the Mersenne prime is 
in the middle are composed by the twins where the 
Mersenne prime is last followed by twins where it 
comes first.   

\subsection{Fermat triplets}

We give a few triplets where the Fermat prime comes 
{\it first.} 

{\it Example~10.}

$2^{2^n}+1, 2^{2^n}+3,~2^{2^n}+9$ yields $3,~5,~11;~
5,~7,~13$ for $n=0,~1.$

$2^{2^n}+1, 2^{2^n}+5,~2^{2^n}+11$ yields $3,~7,~13$ 
for $n=0.$
 
{\bf Corollary~3.11} {\it (i)} $2^{2^n}+1, 2^{2^n}+3,~
2^{2^n}+5$ {\it for} $n=0,$ {\it yields} $3,~5,~7$ 
{\it which is the only case. (ii)} $2^{2^n}+1, 
2^{2^n}+5,~2^{2^n}+9$ {\it for} $n=0$ {\it yields} 
$3,~7,~11$ {\it as the only case.} 
   
{\it Proof.} (i) follows from the factorization (i) 
and (ii) from (ii) in Prop~3.7.~$\diamond$   

We list a few triplets where the Fermat prime comes 
{\it last.} 

{\it Example~11.}

$2^{2^n}-5, 2^{2^n}-3,~2^{2^n}+1$ yields $11,~13,~17$ 
for $n=2.$

$2^{2^n}-9, 2^{2^n}-3,~2^{2^n}+1$ yields $7,~13,~17$ 
for $n=2.$

$2^{2^n}-11, 2^{2^n}-3,~2^{2^n}+1$ yields $5,~13,~17$ 
for $n=2.$

$2^{2^n}-13, 2^{2^n}-3,~2^{2^n}+1$ yields $5,~13,~17$ 
for $n=2.$

Triplets where the Fermat prime is in the middle 
are composed by twins where the Fermat prime is 
last followed by twins where it comes first. 

\section{Regular prime multiplets} 

Extensions of the twin and triplet primes at constant 
distances to quadruplets, quintets, etc. exist but 
are too numerous to be analyzed systematically here.

\subsection{Regular multiplets from quadratic polynomials}

We restrict our attention to those with regularly 
increasing (or decreasing) distances, such as $2N,~N=1, 
2, \ldots$ i.e. $p_1, p_2=p_1+2, p_3=p_1+6, \ldots, 
p_{n+1}=p_1+n(n+1),\ldots, N$ with the sequence of 
differences $\Delta p_j=p_{j+1}-p_j=2j,~ j=1, 2,
\ldots, N.$ Even regular prime triplets and 
quadruplets are too numerous to be listed. We 
therefore start with quintuplets in Theor.~4.1 
below.  

{\it Example~12.} There are at least 14 sextets 
$11, 13, 17, 23, 31, 41;~17, 19,\\23, 29, 37, 47;
~41, 43, 47, 53, 61, 71;~1277, 1279, 1283, 1289, 
1297, 1307;~1607,\\1609, 1613, 1619, 1627, 1637;
~3527, 3529, 3533, 3539, 3547, 3557;~21557, 
21559,\\21563, 21569, 21577, 21587;~28277, 28279, 
28283, 28289, 28297, 28307;\\~31247, 31249, 31253, 
31259, 31267, 31277;~33617, 33619, 33623, 33629, 
33637,\\33647;~55661,\dots,55691;~113147, 113149, 
113153, 113159, 113167, 113177;\\128981,\ldots, 
129011; 548831,\ldots, 548861;~566537,\ldots, 
566567;$\\seven septets $11,\ldots, 41, 
53;~17,\dots, 59;~41,\ldots, 83;~1277,\ldots, 
1319;~21577,\\\dots, 21599;~28277,\dots, 28319;~
113147,\ldots, 113189; 128981,\dots, 129023;$ 
five\\octets $11,\ldots, 53, 67;~17,\ldots, 73;
~41,\ldots, 97;~21557,\ldots, 21599, 21613;\\
128981,\ldots, 129037;$~three nonets $11,\ldots, 
67, 83;~17,\ldots, 89;~41,\ldots, 113$\\and three 
decuplets $11,\ldots, 83, 101;~17,\dots, 117;~41,
\dots, 131.$\\For $N=15$ the sequence $17, 19, 
23, 29, 37, 47, 59, 73, 89, 107, 127, 149, 173,\\
199, 227, 257$ is a regular 16-plet. The second 
such 16-plet starts soon after with 
$41, 43, 47, 53,\ldots, 281=41+15\cdot 16.$ In 
fact, this one extends much longer and is the 
first 40-plet ending with $1523, 1601=41+39
\cdot 40.$ There are also many almost-regular 
prime multiplets where just one member is missing, 
e.g. $n=0.$   

Naturally, questions arise: Are there infinitely 
many of these long regular multiplets or even 
longer ones? Are there other long prime multiplets 
without the regular structure imposed by a (quadratic) 
polynomial on its first multiplet which then 
continues through all of them? Needless to say, 
the existence of so many regular prime 
multiplets linking primes with each other in 
interlocking ways belies the probabilistic 
independence of prime numbers underlying many 
conjectures~\cite{hw},\cite{ms} about them. It may 
not be unreasonable to expect that most multiplets 
repeat infinitely often. Since they are interlocked  
their asymptotic distributions are not independent. 
This suggests that asymptotic laws for prime 
multiplets differ fundamentally from ordinary 
prime numbers. 

The long regular multiplets of Example~12 are 
related to Euler's prime number generating 
polynomials which, in turn, are related to 
imaginary quadratic number fields over the 
rationals. There are corresponding polynomials 
whose values form regular multiplets that are 
related to real quadratic number fields. 
Although some of these polynomials have been 
known for a long time with their large number 
of prime values as the main focus, their regular 
distances within the coherent structure of a 
regular prime number multiplet seem not to have 
been noted (in print). Except for the multiplet 
aspects many details below are known and 
documented in Ref.~\cite{pr}, but some are new.       
  
{\bf Theorem~4.1.} {\it (i) The Euler polynomials 
$E_p(x)=x^2+x+p$ with the prime numbers $p=2, 3, 5, 
11, 17, 41$ assume prime number values $p+x(x+1)$ 
at distances $2(x+1)$ for $x=0, 1, \ldots p-2$ 
forming a regular $(p-1)-$plet. 

(ii) The polynomials $f_p(x)=2x^2+p$ with the 
primes $p=3, 5, 11, 29$ assume prime number 
values at $x=0, 1, \dots, p-1$ forming a regular 
$p-$plet at distances $2(2x+1).$ 

(iii) The polynomials $F_p(x)=2x^2+2x+\frac{1}{2}
(p+1)$ with primes $p=5, 13, 37\equiv 1\pmod{4}$ 
assume prime values for $x=0, 1, (p-3)/2$ which 
form a regular $(p-1)/2-$plet at distances $4(x+1).$ 

(iv) The polynomials $G_p$, with primes $p, q$, 
\begin{eqnarray}\nonumber
G_p(x)&=&px^2+px+\frac{1}{4}(p+q),~p<q, pq\equiv 3 
\pmod{4},\\\nonumber (p,q)&=&(3,5), (3,17), (3,41), 
(3,89), (5, 7), (5, 23), (5, 47), (7, 13), (7, 61),\\&&
(11, 17), (13,31)  
\end{eqnarray} 
assume prime values for $x=0, 1, \ldots \frac{1}{4}
(p+q)-2.$ which form regular multiplets at distances 
$2p(x+1)$ that are independent of the prime $q.$ 

(v) Polynomials corresponding to real quadratic number 
fields are $g_d(x)=-x^2+x+\frac{1}{4}(d-1)$ with 
$d>0, d\equiv 1 \pmod{4}$ and square-free. Their 
values are prime numbers for $x=2, 3, \ldots <
\frac{1}{2}\sqrt{d-1}$ forming regular multiplets 
at distances $2p(x+1).$ Relevant values are $d=37, 
53, 77, 101, 173, 197, 293, 437,\\677.$ 

(vi) Quintets, sixtets, septets, octets, nonets and 
decuplets are generated by the following polynomials}  
\begin{eqnarray}\nonumber
Q_p(x)&=&x^2+x+p,~x=0,\ldots,4;~p=11, 17, 41, 347, 641, 
1427, 4001,\\\nonumber&&4637, 4931, 19421, 22271, 
23471, 26711, 27941, 28277,\\\nonumber&&31247, 32057, 
33617; 113147\ldots;\\\nonumber SX_p(x)&=&x^2+x+p,~x=0,
\ldots 5;~p=11, 17, 41, 1277, 1607, 3527,\\
\nonumber&&28277, 31247, 33617, 55661, 113147, 128981, 
548831, 566537,\ldots;\\\nonumber~S_p(x)&=&x^2+x+p,
~x=0,\ldots, 6;~p=11, 17, 41, 1277, 28277,\\
\nonumber &&113147, 128981,\ldots;\\\nonumber~O_p(x)&=&
x^2+x+p,~x=0,\ldots, 7;~p=11, 17, 41, 128981,
\ldots;\\\nonumber~N_p(x)&=&x^2+x+p,~x=0,\ldots, 8;~p=
11, 17, 41,\ldots;\\~D_p(x)&=&x^2+x+p,~x=0,\ldots, 
9;~p=11, 17, 41,\ldots .
\end{eqnarray}

{\it Proof.} (i) It is long known~\cite{pr} that for 
the listed primes $p$ the values $E_p(x)=x^2+x+p$ at 
$x=0, 1, 2, \ldots, p-2$ are prime numbers. Since 
\begin{eqnarray}
E_p(x+1)-E_p(x)=(x+1)^2-x^2+1=2(x+1),
\end{eqnarray}   
the primes $p+x(x+1)$ form a regular $(p-1)-$plet 
at distances $2(x+1)$ from each other for $x=0, 1, 
2, \ldots, p-2.$ 

(ii) Since $f_p(x+1)-f_p(x)=2(2x+1),$ the primes 
$p+2x^2,~x=0, 1, \ldots, p-1$ form a regular 
$p-$plet at distances $2(2x+1).$ 

(iii) Since $F_p(x+1)-F_p(x)=4(x+1)$ the primes 
$2x(x+1)+(p+1)/2$ form a regular $(p-1)/2-$plet 
at distances $4(x+1).$ 

(iv) Since $G_p(x+1)-G_p(x)=2p(x+1)$ the multiplet 
structure is clear.   

(v) For $g_d,$ this follows from $g_d(x+1)-g_d(x)
=-2x.$ These are multiplets with linearly 
decreasing distances. 

(vi) These multiplets may be verified by a table of 
prime numbers or symbolic-math software.$~\diamond$ 
    
There are many more recently found polynomials in 
the literature~\cite{pr} which also form regular 
prime multiplets.   
 
\subsection{Optimal quadratic polynomials}

This subject has a long history~\cite{dick} with 
a rather unsystematic record. It is well known that, 
if $P(x)=\sum_{j=0}^n a_n x^n$ is a non-constant 
polynomial of degree $n\geq 1$ with integral 
coefficients $a_j,~|a_0|=p_0$ a prime number, then 
$P(x)$ can assume prime values at most for $x=0, 
1,\ldots, p_0-1$ because $p_0|P(p_0)$.   

{\bf Definition~4.2.} The polynomial $P(x)=
\sum_{j=0}^n a_j x^j$ is called {\it optimal} if 
$|P(j)|=p_j$ is prime for $j=0, 1,\dots, p_0-1,$ 
forming a regular $p_0-$plet. 

Prime number values have to be successive, but 
they may repeat and be negative. This is often 
caused by negative coefficients in a polynomial.  
 
As $E_p(p-1)=p^2,$ the Euler polynomials in (i) 
of Theor.~4.1 for $p=2, 3, 5, 11, 17, 41$ are 
one prime value short of being optimal. But  
Legendre's quadratic polynomials~\cite{dick} 
$f_p(x)=2x^2+p$ for $p=3, 5, 11, 29$ in (ii) 
of Theor.~4.1 are optimal.     

Since $E_p(-x)=E_p(x-1),$ Euler polynomials 
give repeating prime multiplets when they are 
considered for positive and negative argument. 
For more general polynomials this is not 
the case. Upon shifting the argument, the 
modified Euler polynomials $E_p(x-1)=x(x-1)+p$ 
do become optimal, but they repeat the initial 
prime value. More generally, the identity
\begin{eqnarray}
E_q(x-n)=x^2+(1-2n)x+E_q(n-1)
\end{eqnarray}  
for prime numbers $q,~p=E_q(n-1)$ and nonnegative 
integer $n$ leads to many new repeating polynomials 
with more numerous prime values. For $q=41$ and 
$n=2, 3,\ldots, 40$ none of these polynomials is 
optimal (cf. Eq.~(\ref{43})) in Cor.~4.5 e.g.), 
though, including Escott's for $n=40$~\cite{pr},
\cite{dick}. 
 
{\bf Definition~4.3.} The polynomial $P(x)$ is 
{\it bi-optimal} if $|P(j)|$ are prime for 
$1-p_0\leq j\leq p_0-1,$ forming a $(2p_0-1)-$plet.

This corresponds to combining $P(j)$ for $j=0, 1,
\ldots, p_0-1$ and $P(-j)$ for $j=1,\ldots, p_0-1$ 
into one prime number multiplet. For 
polynomials of odd degree the starting value 
$x=0$ is appropriate. For polynomials of even 
degree it is perhaps more natural to include 
negative arguments as well. 

We now address a question raised by Legendre's 
$f_p(x)$ and Euler's modified polynomials: Are 
there optimal quadratic polynomials for the 
missing-link primes $p=7, 13, 19, 23, 31, 37?$

For $p=7,~Q_2(x)=2x(x-1)+7$ is optimal, if 
repeating, generating a regular prime number 
septet with distances $\Delta Q_2(x)=4x,$ and 
it is almost bi-optimal because $Q_2(-x)=Q_2(x+1)$. 

For $p=19,~Q_7(x)=2x(x-1)+19$~\cite{dick} is 
not only optimal but almost bi-optimal in view 
of $Q_7(-x)=Q_7(x+1)$, generating a regular 
36-plet with the same distance law as $Q_2.$ 

For $p=23,~Q_8(x)=3x(x-1)+23$~\cite{dick} is 
not only optimal but almost bi-optimal forming 
a 42-plet because $Q_8(-x)=Q_8(x+1)$.   
     
{\bf Proposition~4.4.} {\it (i) There is an 
infinity of optimal quadratic polynomials for} 
$p_0=2:$ 
\begin{eqnarray}
Q_1(x)=x^2+(p_1-3)x+2,
\end{eqnarray} 
{\it where $p_1$ is prime. If $|p_1-6|$ is 
prime then $Q_1(x)$ is bi-optimal.   
 
(ii) There are at least three quadratic 
polynomials for} $p_0=13:$
\begin{eqnarray}\nonumber
Q_3(x)&=&x^2+27x+13,~Q_4(x)=x^2-3x+13,\\
Q_5(x)&=&2x^2-4x+13,
\end{eqnarray}
{\it that form 12-plets, one prime value short 
of optimal. Hence  
\begin{eqnarray}
Q_3(x-1)=x^2+25x-13
\end{eqnarray}
is optimal. $Q_5(-x)$ forms a decuplet making 
$Q_5(x)$ almost bi-optimal. There is another 
polynomial,  
\begin{eqnarray}
Q_6(x)=2x^2+26x+13, 
\end{eqnarray}
that forms a decuplet.}

{\it Proof.} (i) follows from $Q_1(0)=2,~
Q_1(1)=p_1$ and $Q_1(-1)=6-p_1.$    

(ii) For monic optimal polynomials let 
\begin{eqnarray}
Q(x)=x^2+bx+p_0,~Q(1)=p_1=p_0+b+1,~Q(2)
=p_2=p_0+4+2b,
\label{qm1}
\end{eqnarray}
where $p_0, p_1, p_2$ are prime numbers. Then 
\begin{eqnarray}
b=p_1-p_0-1,~Q=x^2+(p_1-p_0-1)x+p_0,~p_2
=2p_1-p_0+2.
\label{qm2}
\end{eqnarray}
If $Q(x)=2x^2+bx+p_0,$ then 
\begin{eqnarray}
Q(1)=p_1=p_0+b+2,~Q(2)=p_2=p_0+2b+8. 
\label{q2}
\end{eqnarray}

(ii) For $p_0=13$ and $p_1=41$ we obtain $Q_3$ 
from Eq.~(\ref{qm2}) and the 12-plet  
\begin{eqnarray}\nonumber
Q_3(1)&=&41,~Q_3(2)=71,~Q_3(3)=103,~Q_3(4)=137,~
Q_3(5)=173,\\\nonumber Q_3(6)&=&211,~Q_3(7)=251,~
Q_3(8)=293,~Q_3(9)=337,~Q_3(10)=383,\\Q_3(11)&=&431, 
\end{eqnarray}
one prime value short of optimal. $Q_3$ is ascending 
and non-repeating. 

For $p_0=13,~p_1=11$ we get the monic $Q_4,$ again 
with a 12-plet of prime values 
\begin{eqnarray}\nonumber
Q_4(1)&=&11,~Q_4(2)=11,~Q_4(3)=13,~Q_4(4)=17,~Q_4(5)
=23,\\\nonumber Q_4(6)&=&31,~Q_4(7)=41,~Q_4(8)=53,~
Q_4(9)=67,~Q_4(10)=83,\\Q_4(11)&=&101,
\end{eqnarray}
which is one prime value short of optimal. It is 
repeating and non-ascending because $Q_4(1)=11=
Q_4(2),~Q_4(3)=13=Q_4(0).$

For $p_1=11$ we get $Q_5$ from Eq.~(\ref{q2}) 
and the 12-plet 
\begin{eqnarray}\nonumber
Q_5(1)&=&11,~Q_5(2)=13,~Q_5(3)=19,~Q_5(4)=29,~Q_5(5)
=43,\\\nonumber Q_5(6)&=&61,~Q_5(7)=83,~Q_5(8)=109,~
Q_5(9)=139,~Q_5(10)=173,\\Q_5(11)&=&211,
\end{eqnarray}
which is one prime value short of optimal. It is 
obviously repeating and non-ascending. Since 
\begin{eqnarray}
Q_5(-x)=2x^2+4x+13=Q_5(x+2), 
\end{eqnarray}
it generates a decuplet. Therefore, $Q_5(x)$ is 
almost bi-optimal.  

For $p_1=41$ we get $Q_6$ from Eq.~(\ref{q2}) and 
the decuplet 
\begin{eqnarray}\nonumber
Q_6(2)&=&73,~Q_6(3)=109,~Q_6(4)=149,~Q_6(5)=193,~
Q_6(6)=241,\\Q_6(7)&=&293,~Q_6(8)=349,~Q_6(9)=409.
~\diamond
\end{eqnarray} 

{\bf Corollary~4.5.} {\it The polynomial} 
\begin{eqnarray}
Q_9(x)=x^2+3x+19  
\end{eqnarray}
{\it forms a 15-plet for} $x=0,\ldots, 14$,    
\begin{eqnarray}
Q_{10}(x)=x^2-x+11,~Q_{11}(x)=2x^2+22x-11
\end{eqnarray}
{\it are optimal polynomials for} $11,$ {\it 
and}
\begin{eqnarray}
Q_{12}(x)=2x^2-4x+31,~Q_{13}(x)=x^2-3x+43
\label{43}
\end{eqnarray}
{\it form 30- and 42-plets, respectively, one 
prime value short of optimal.}

{\it Proof.} This follows from Prop.~4.4 in 
conjunction with  
\begin{eqnarray}\nonumber
Q_9(x-1)&=&x^2+x+17=E_{17}(x),~Q_{10}(x)=E_{11}(x-1),
\\\nonumber Q_{11}(x)&=&Q_6(x-1),~Q_{12}(x)
=f_{29}(x-1),\\Q_{13}(x)&=&E_{41}(x-2), 
\end{eqnarray} 
respectively.~$\diamond$

We leave open the questions: Are there primes for 
which there is no optimal quadratic polynomial?   
Are $p=31,~37$ such cases? 

\subsection{Optimal cubic polynomials}

We now investigate cubic polynomials for the prime 
numbers $p_0=2, 3, 5, 7, 11,\\13$. 

{\bf Theorem~4.6.} {\it (i) There is an infinity 
of optimal cubic polynomials for the prime numbers 
$2, 3$. For} $p_0=2:$ 
\begin{eqnarray}
C_1(x)=x^3+mx^2+(p_1-3-m)x+2,~C_1(1)=p_1,
\end{eqnarray} 
{\it where $m$ is an arbitrary integer and $p_1$ 
a prime number. If $|4+2m-p_1|$ is prime then 
$C_1(x)$ is bi-optimal. For} $p_0=3:$
\begin{eqnarray}\nonumber
C_2(x)&=&x^3+[\frac{1}{2}(p_2-3)-p_1]x^2+[2p_1-
\frac{1}{2}(p_2+5)]x+3,\\C_2(1)&=&p_1,~C_2(2)=p_2,
\end{eqnarray}
{\it where $p_1, p_2$ are prime and $p_2$ odd. If 
$|p_2-3p_1+3|, |3p_2-8p_1-6|$ are prime then $C_2(x)$ 
is bi-optimal. 

(ii) There are optimal polynomials for} $p_0=5:$ 
\begin{eqnarray}
C_3(x)=x^3-x^2+2x+5,~C_4(x)=2x^3+4x^2-4x+5.
\end{eqnarray}
$p_0=7:$
\begin{eqnarray}
C_5(x)=x^3-x^2+6x+7,
\end{eqnarray}
$p_0=11:$
\begin{eqnarray}
C_6(x)=x^3+5x^2+2x+11,~C_7(x)=x^3-4x^2+5x+11,
\end{eqnarray}
$p_0=13:$
\begin{eqnarray}
C_8(x)=x^3-5x^2+8x+13.
\end{eqnarray}
{\it Proof.} (i) $C_1(1)=p_1$ is readily 
verified for $p_0=2$ and $C_2(1)=p_1,~C_2(2)
=p_2$ for $p_0=3.$ Note that $C_1(-1)=4+2m-p_1, 
C_2(-1)=p_2-3p_1+3, C_2(-2)=3p_2-8p_1-6.$      

(ii) Let $C(x)=ax^3+bx^2+cx+p_0.$ Then 
\begin{eqnarray}\nonumber
p_1&=&a+b+c+p_0,~p_2=8a+4b+2c+p_0,\\p_3
&=&27a+9b+3c+p_0.
\end{eqnarray}
(ii) For $p_0=5,$ choosing $a=1, p_1=7, p_2=13$ 
yields $b+c=1,~2b+c=0$ and therefore $C_3(x)$ 
yielding the following additional prime values 
\begin{eqnarray}
C_3(3)=29,~C_3(4)=61,
\end{eqnarray} 
forming altogether a quintet. 

For $a=2,p_0=5,~p_1=7,~p_2=29$ we get $c=-b$ 
and $b=4.$ The polynomial $C_4$ yields the 
additional prime values
\begin{eqnarray}
C_4(3)=83,~C_4(4)=181,  
\end{eqnarray} 
forming a quintuplet. For $p_0=7,~p_1=13,
~p_2=23$ we get $b=-1,~c=6$ and $C_5$ yields 
the additional prime values 
\begin{eqnarray}
C_5(3)=43,~C_5(4)=79,~C_5(5)=137,~C_5(6)=223, 
\end{eqnarray} 
a septet and optimal. For $p_0=11,~p_1=19,~p_2=43$ 
we get $b=5,~c=2$ and $C_6$ yields the prime values 
\begin{eqnarray}\nonumber
C_6(3)&=&89,~C_6(4)=163,~C_6(5)=311,~C_6(6)=419,~
C_6(7)=613,\\C_6(8)&=&859,~C_6(9)=1163,~C_6(10)=1531, 
\end{eqnarray}
forming an optimal 11-plet. It may be verified 
that $C_7$ forms an 11-plet also. $C_8(x)$  
generates the optimal 13-plet 
\begin{eqnarray}\nonumber
C_8(0)&=&13,~C_8(1)=17,~C_8(2)=17,~C_8(3)=19,~
C_8(4)=29,\\\nonumber~C_8(5)&=&53,~C_8(6)=97,~
C_8(7)=167,~C_8(8)=269,~C_8(9)=409,\\C_8(10)&=&
593,~C_8(11)=827,~C_8(12)=1117.~\diamond
\end{eqnarray}

Concluding we ask: Are there optimal cubic 
polynomials for all prime numbers $p>13?$


\begin{thebibliography}{0}  

\bibitem{hjw} Weber, H. J., 2010, ``Generalized 
Twin Prime Formulas,'' Global J. of Pure and Applied 
Math. 6(1), pp. 101-116.  

\bibitem{pr} Ribenboim, P., 1996, The New Book 
of Prime Number Records, Springer, Berlin , and refs. 
therein. 

\bibitem{hw} Hardy, G. H. and Wright, E. M., 1988, 
An Introduction to the Theory of Numbers, 
Clarendon Press, Oxford, 5th ed. 

\bibitem{ms} Schroeder, M. R., 1997, Number 
Theory in Science and Communication, 3rd ed., 
Springer, Berlin. 

\bibitem{dick} Dickson, L. E., 2005, A History of 
Number Theory, Dover, New York, Vol. 1, pp. 420-421. 

\end{thebibliography}
\end{document}